# LIMITED RESOURCE OPTIMAL DISTRIBUTION ALGORITHM BASED ON GAME ITERATION METHOD

V.Ya. Vilisov
Department of Mathematics
University of Technology
Korolyov city, Russia
vvib@yandex.ru

*Abstract* — *The article provides a solution algorithm for the linear programming problem (LPP) with the latter being presented as an antagonistic matrix game so the game's further solution is based on the iterative method. The algorithm is presented as a computer program. Having applied necessary accuracy, the author has researched the solution assessment convergence rate in relation to the actual value. Program implementation demonstrates high rate of the LPP solution receipt, with the acceptable accuracy being fractions or unities. It allows using the algorithm in the integrated systems for the purpose of their optimal control.*

## INTRODUCTION

The task of limited resource distribution can be met rather often in different applied disciplines. Usually, a choice of optimal control in the modern control systems occurs with a certain regularity or as a reaction at some events. A majority of the situations involving the choice of optimal control can be presented as mathematical programming problems, with the linear programming problems (LPP) amounting to their considerable proportion.

LPP takes a special place among the methods of optimal control and/or optimal distribution of limited resources as it can become a conclusion of many applied problems whose original form does not resemble LPP at all; they are, for example, controlled Markov chains [1], transport problem [2, 3], problem of allocation [3], matrix games [4, 5] etc. The applied industries, which are used to solve the problems, are also very diverse. They include, for instance, control of robotic groups [6], distribution of communication channels in robotic systems (RS), distribution of efforts and means when managing emergency situations, building an optimal industrial program [7] etc.

John von Neumann, a matrix game founder, identified [4] a close connection between matrix games and linear programming [8]. However, in practice, we are mostly aware of those methods of solving antagonistic matrix games that are reduced to LPP, thus allowing finding the solution of the random dimensionality matrix games.

Simplex method is a standard way of solving LPP (in various modifications). Still, multiple research [8] demonstrates that it cannot always provide a guarantee of finding an optimal solution. Other research methods also possess this defect. In many ways, the convergence indicators of research methods depend on the data of a specific applied problem. Heuristic modifications of simplex method, for example those applied in MS Excel, are usually directed at the convergence improvement. Moreover, other disadvantages of this method include awkwardness and non-obviousness of computational procedures. Today, these disadvantages are usually compensated with the user qualification, who, when working with the majority of modern optimization software, is presented with a range of parameters that can be modified for the purpose of obtaining an acceptable solution. However, all attempts that have been made to implement the above-mentioned models in the real-time control systems, for instance in on-board computers of modern integrated systems, have led to a significant decrease of their operating reliability. Therefore, we have a question: can we find a replacement to our "workhorse", i.e. the simplex method?

The algorithm, which has been suggested in this article, is a considerably more simple and obvious method that allows to obtain a solution with a prescribed accuracy. This provides a possibility, in each specific application, to find a compromise between the solution's accuracy and number of iterations, i.e. the computational speed. The method offered herein is based on the one-one correlation of LPP and antagonistic matrix games (AMG) and bears the name of the **Game Iteration Method**. Let us consider its main elements.

## LPP-BASED TRADITIONAL AMG SOLUTION SCHEME

AMG (zero sum games) have been traditionally [5] solved by their presentation as LPP. The essence of this transformation is the following. Usually, the game is set with a payoff matrix $\tilde{A} = [\tilde{a}_{ij}]_{mn}$. $\tilde{a}_{ij}$ elements demonstrate the *first player*'s payoff (P1) in each game. Since the game has a zero sum, $\tilde{a}_{ij}$ is simultaneously also the loss value of the *second player* (P2). *To solve the game* means [5] to find a pair of *optimal mixed strategies* of players $\bar{p} = [p_1 \; p_2 \; ... \; p_m]^T, \bar{q} = [q_1 \; q_2 \; ... \; q_n]^T$, as well as the *game's price* $V$. Here, $\bar{p}$ and $\bar{q}$ are vector-columns of dimensionality, correspondingly of $m$ and $n$, while $T$ is a conjugation symbol. From now on we shall use a *tilde* to mark the elements of AMG and LPP that belong to the game form (solution of AMG with the aid of LPP). Using the vector-matrix form, the game price can be presented as $V = \bar{p}^T \tilde{A} \bar{q}$. The normalization condition shall be performed for the vector components $\bar{p}$ and $\bar{q}$ that include the likelihood of players may applying some pure strategy:

$$\sum_{i=1}^{m} p_i = 1; \qquad (1)$$

$$\sum_{j=1}^{n} q_j = 1. \tag{2}$$

If P1 applies his optimal mixed strategy and P2 any pure one, P1 shall receive the payoff not less than the $V$ game price should the game be repeated multiple number of times:

$$\sum_{i=1}^{m} \tilde{a}_{ij} p_i \geq V; \; j = \overline{1,n}. \tag{3}$$

Having introduced the symbol: $\tilde{x}_i \triangleq \frac{p_i}{V}$ and having inserted it in (1) we shall receive $\sum_{i=1}^{m} \tilde{x}_i = \frac{1}{V}$. The P1's objective is to maximize his average payoff, i.e. $V$, the game price. It means that he should minimize the $\frac{1}{V}$ value. The values without a tilde shall mean variables of a standard LPP and corresponding indices of the delimitation system that are not directly related to the AMG. Therefore, with the consideration of introduced symbols and after having replaced variables in (3), we can write down the AMG-obtained LPP setting as follows:

$$\sum_{i=1}^{m} \tilde{a}_{ij} \tilde{x}_i \geq 1; \; j = \overline{1,n}, \tag{4}$$

$$L(\bar{\tilde{x}}) = \sum_{i=1}^{m} \tilde{x}_i \to \min_{\tilde{x}_i}, \tag{5}$$

where $L(\bar{\tilde{x}})$ is a LPP's objective function based on AMG.

The solution of the LPP shall lead to the following solution of AMG for P1:

$$p_i = \frac{\tilde{x}_i}{L_{opt}}; \; i = \overline{1,m}; \; V = \frac{1}{L_{opt}}, \tag{6}$$

where $L_{opt}$ is an optimal value of the objective function and $\tilde{x}_i$ are the elements of the P1's optimal mixed strategy vector.

We can also perform similar constructions for the second player (P2), with the result being the following:

$$\sum_{j=1}^{n} \tilde{a}_{ij} \tilde{y}_j \leq 1; \; i = \overline{1,m}, \tag{7}$$

$$L(\bar{\tilde{y}}) = \sum_{j=1}^{n} \tilde{y}_j \to \max_{\tilde{y}_j} \tag{8}$$

The problem possesses dual characteristics in relation to the source LPP (4)-(5). After solving it, we shall obtain an optimal game solution for P2:

$$q_j = \frac{\tilde{y}_j}{L_{opt}}; \; j = \overline{1,n}; \; V = \frac{1}{L_{opt}}, \tag{9}$$

That is a traditional technology [5, 9, 10] of solving AMG by its conversion into LPP, its solution and ultimate transformation of the obtained LPP solution into the AMG solution. At that, it is assumed that *there exists an acceptable method of solving LPP*, which is usually a simplex method [3]. Next we shall demonstrate a reverse algorithm, i.e. determination of the LPP solution by its transformation into AMG, searching for the solution and performing a reverse transformation of an obtained AMG solution into the LPP solution.

## AMG-BASED LPP SOLUTION

In order to apply a reverse algorithm, we shall possess a matrix game solution method that is *not based on LPP*. There are not so many methods like that: method of Brown-Robinson [9, 10], monotonic iterative algorithm [11], graphical solving method [3] and some others. The most universal is a Brown-Robinson method, which allows looking for the solution of a random dimensionality AMG. Its structure is rather simple, but it can require a large number of monotypic calculations should the dimensionality of an AMG payoff matrix be very big. As the computing powers of computers (including on-board ones) are constantly growing, this method is becoming rather attractive for the AMG solution.

Here we shall demonstrate the calculations that allow providing a LPP solution based on the solution of an equivalent AMG. At that, it is assumed that the AMG solution has been obtained without using the LPP solution algorithms. This scheme allows for the immediate receipt of both direct and dual LPP.

Let us presume that LPP that has to be solved is presented as a maximum problem (thus the dual problem that corresponds to it shall be presented as a minimum problem). There are several forms of presenting LPP [8], which can be transformed into each other using a one-one principle. Next we shall apply a form of the maximum problem with the delimitations being ($\leq$):

$$L(\bar{y}) = \bar{c}^T \bar{y} = \sum_{j=1}^{n} c_j y_j \to \max_{y_j}, \tag{10}$$

$$\sum_{j=1}^{n} a_{ij} y_j \leq a_{i0}; \; i = \overline{1,m}; \; y_j \geq 0, \tag{11}$$

where $\bar{y} = [y_1 \; y_2 \; \ldots \; y_n]^T$ is a vector of variables;

$\bar{c} = [c_1 \; c_2 \; \ldots \; c_n]^T$ - is an index vector of a linear objective function;

$A = [a_{ij}]_{mn}$ - is an index matrix of the linear delimitation system;

$\bar{a}_0 = [a_{10} \; a_{20} \; \ldots \; a_{m0}]^T$ - is a constant term vector of the delimitation system.

Let us presume that $c_j > 0, a_{ij} > 0, a_{i0} > 0$. It is easy to see that the problem (10)-(11) can be reduced to the problem (7)-(8) should we introduce the following symbols:

$$\tilde{a}_{ij} = \frac{a_{ij}}{a_{i0} c_j}; \tag{12}$$

$$\tilde{y}_j = c_j y_j \text{ so that } y_j = \frac{\tilde{y}_j}{c_j}. \tag{13}$$

There is a one-one correlation between the (7)-(8) problem and the payoff matrix $\tilde{A} = [\tilde{a}_{ij}]_{mn}$, whose elements are fully defined by the correlation (12). Therefore, it is possible to conclude that AMG was built on the basis of the LPP source data.

In order to solve AMG, whose payoff matrix consists of elements mentioned in (12), we shall use a method of Brown-Robinson [9, 10], with the result being two vectors of optimal mixed strategies and the game price: $\bar{p} = [p_1 \; p_2 \; \ldots \; p_m]^T, \bar{q} = [q_1 \; q_2 \; \ldots \; q_n]^T, V$.

Having obtained an AMG solution and using (6), (9) and (13) it is possible to obtain a solution of the direct LPP, dual LPP as well as the optimal value of the objective function pertaining to the source LPP (they are similar both for the direct and dual problems):

$$L_{opt} = \frac{1}{V}; \quad y_j^{opt} = \frac{q_j}{c_j V}; \quad x_i^{opt} = \frac{p_i}{a_{io} V}, \quad (14)$$

where $i = \overline{1,m}; \; j = \overline{1,n}$.

Thus, sequence of steps leading to the LPP solution (10) - (11) is the following:

1. Construction of the AMG payoff matrix on the basis of the source LPP's parameters via the calculation of the AMG elements with the formula (12).

2. Solving AMG using a Brown-Robinson method.

3. Recalculation of an obtained AMG solution into the LPP solution using the formulas (14).

Let us note some peculiarities of transforming LPP into AMG and back:

*a)*. if LPP is a *maximum* with *n* variables (columns) and *m* delimitations (lines), the index matrix of the left-side delimitations should be transformed into a payoff matrix that would also possess *m* lines and *n* columns. After solving the game, optimal values of *n* LPP variables can be obtained from the *n*-dimensional vector related to the optimal mixed strategy of the *second player*;

*b)*. if LPP is a *minimum* with *n* variables (columns) and *m* delimitations (lines), the index matrix of the left-side delimitations should be transformed into a payoff matrix and conjugated so that it would also possess *n* lines and *m* columns. After solving the game, optimal values of *n* LPP variables can be obtained from the *n*-dimensional vector related to the optimal mixed strategy of the *first player*.

## METHOD OF BROWN-ROBINSON

Whereas in the method of game iterations applied for the AMG solution we have used Brown-Robinson method, let us provide here its summary. It was suggested by George W. Brown [9], while its convergence to the optimal (minimax) solution was proved by Julia Robinson [10]. It is also known as a *method of fictitious play* or *iterative method*. It is an iterative procedure that imitates a game looping, while it is also assumed that players perform an alternate choice of their next pure strategy taking into account all information on actions of an opponent that they have made in the previous games. Estimations of optimal mixed strategies are computed in the form of current frequencies related to the players' usage of their pure strategies during the whole observation interval, while the estimation of the game price is computed as an average of a current payoff of the P1 player and a loss of the P2 player. The iterations stop when the set estimation accuracy has been achieved. Let us provide a formal presentation of this procedure.

Let us assume that two players know all elements of the payoff matrix $\tilde{A} = [\tilde{a}_{ij}]_{mn}$. They play several games, provided that in each game both players are able to see pure strategies that have been chosen by the opponent. In each game players choose their best strategies (P1 to maximize his payoff, P2 to minimize his loss). Performing his next (*t*+1) choice, each player takes into account pure strategy statistics chosen by his opponent during previous *t* steps: if P2 has chosen his *j* strategy $Q_j$ times out of total number *t* of games,

P1 shall choose his *i* pure strategy so that he would maximize an average value of the payoff:

$$i(t+1) = arg \max_{i=\overline{1,m}} v_i(t) = arg \max_{i=\overline{1,m}} \sum_{j=1}^{n} a_{ij} \hat{q}_{j(t)}, \quad (15)$$

where $\hat{q}_{j(t)} \triangleq \frac{Q_{j(t)}}{t}$ is a current (after *t* steps) estimation of an optimal mixed strategy of P2; $v_i(t)$ is an estimation of a current average payoff of P1 should he apply his *i* pure strategy.

Similarly, if P1 has applied his *i* strategy $P_j$ times, P2 shall choose his *j* strategy so that he would minimize an average loss value:

$$j(t+1) = arg \min_{j=\overline{1,n}} v_j(t) = arg \min_{j=\overline{1,n}} \sum_{i=1}^{m} a_{ij} \hat{p}_{i(t)}, \quad (16)$$

where $\hat{p}_{i(t)} \triangleq \frac{P_{i(t)}}{t}$ is a current (after *t* steps) estimation of an optimal mixed strategy of P1; $v_j(t)$ is an estimation of a current average payoff of P2 should he apply his *j* pure strategy.

Random values $\hat{p}_i(t)$ and $\hat{q}_j(t)$ demonstrate sequences of players' mixed strategy vector element estimations that would converge in the limit (as it was proved by J. Robinson [10]) in relation to the optimal mixed strategies. And the estimation of an average payoff of P1 and average loss of P2 that have been computed with the consideration of current estimations of optimal mixed strategies, results in the game price:

$$V = \lim_{t \to \infty} \max_{i=\overline{1,m}} v_i(t) = \lim_{t \to \infty} \min_{j=\overline{1,n}} v_j(t). \quad (17)$$

Iterative process can be stopped when it has achieved necessary accuracy, with the latter being controlled as *an absolute moving average value related to the game price estimation difference* both for P1 and P2:

$$\Delta v(t) = \left| \max_{i=\overline{1,m}} v_i(t) - \min_{j=\overline{1,n}} v_j(t) \right|. \quad (18)$$

When stopping an iterative process we can also use *moving variability* or *moving root-mean-square deviation* of $\Delta v(t)$ value. Moreover, it is possible to build a stopping rule by introducing a specific criterion for the estimation vectors related to the players' optimal strategies. We shall demonstrate it further with an example.

Multiple imitation experiments showed that the convergence process related to the acceptable number of iterations depends strongly on the game's payoff matrix element values, in particular, on the difference between the upper and lower price of the game. Here we demonstrate an estimation convergence for the $3 \times 3$ game.

***Example.*** The payoff matrix is as follows:

$$A = \begin{bmatrix} 6 & 1 & 4 \\ 2 & 4 & 2 \\ 4 & 3 & 5 \end{bmatrix}.$$

The game does not have a saddle point; therefore, mixed strategies should be applied when searching for the solution. If we use a traditional method of game convergence to LPP, the exact solution in this case shall be the following: $\bar{p} = [0 \quad 0.33 \quad 0.67]^T$, $\bar{q} = [0.33 \quad 0.67 \quad 0]^T$, $V = 3.33$.

Graphical demonstration of the convergence for the solution's estimation elements obtained due to the method of Brown-Robinson for 100 steps is shown on Picture 1.

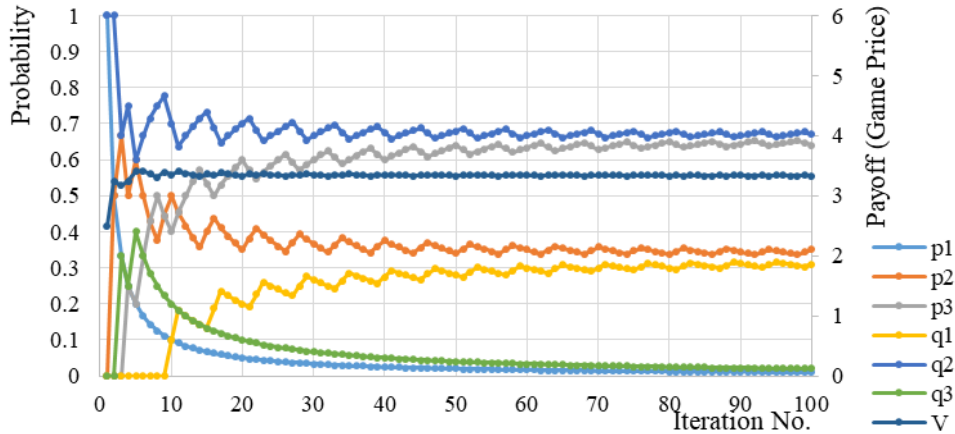

Pic. 1. Estimation convergence of the matrix game solution

These are the final estimation values of the solution per 100 iterations: $\hat{p} = [0.01 \quad 0.35 \quad 0.64]^T$, $\hat{q} = [0.31 \quad 0.67 \quad 0.02]^T$, $\hat{V} = 3.34$. Thus, the deviation of the game price from the exact value is 0.3%.

When dealing with model examples (when we know the exact solution), it is possible to consider a *normalized distance of a vector showing the estimate vector difference and an actual vector of the optimal mixed strategy* of the player as the precision measure of the current estimates related to the elements of the solution. Thus, the vector distance $\hat{\rho}_k = \bar{p} - \hat{p}_k$, being the distance from the normalized (located within the interval $[0; 1]$) estimate vector to the normalized vector of an actual optimal mixed strategy is shown on Picture 2.

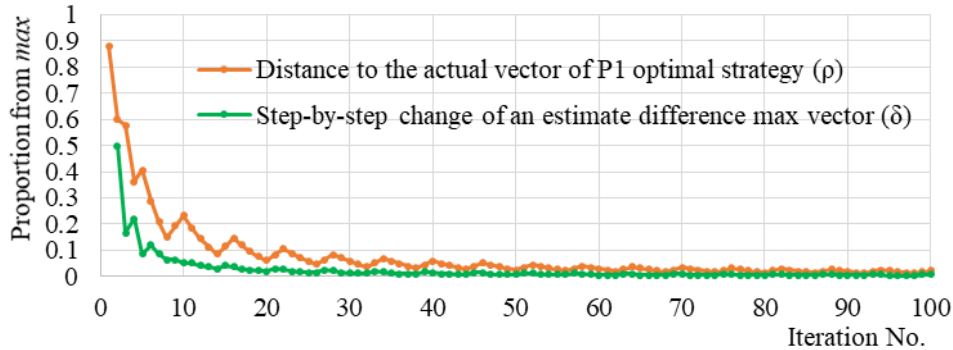

Pic. 2. Convergence of indices pertaining to the solution exactness

The picture also contains another measure, i.e. one computed at each iteration. It is a step-by-step change of the maximum distance (norm) related to the vector showing the estimate vector difference in relation to the optimal strategies of P1 and P2 at adjacent iteration steps: $\hat{\delta}_k = \max(\|\hat{p}_k - \hat{p}_{k-1}\|; \|\hat{q}_k - \hat{q}_{k-1}\|)$. Estimates for the example state that after 100 iterations we have $\hat{\rho}_{100} = 0.023$. It comprises 2.3% off the maximum possible deviation of the estimate vector from the actual one. Still, this measure is available only within the *simulation research conditions*. Another measure of convergence ($\hat{\delta}_k$) does not require to know the value of an actual optimal mixed strategy vector, being, therefore, *available within conditions of the real application of the algorithm*. Its value for the given example is $\hat{\delta}_{100} = 0.007$, i.e. 0.7% off the maximum possible value. In the given example, the correlation coefficient between these two measures constituted 0.87, thus showing a high statistic connection between these two measures of estimate convergence. Therefore, when using the Brown-Robinson algorithm in practice for the purposes of stopping the iterative algorithm, the application of $\hat{\delta}$ index is a highly likely one.

## LPP SOLUTION ALGORITHM USING A GAME ITERATION METHOD

A highly detailed solution of the LPP by its transformation into the AMG, solving AMG using the Brown-Robinson method, as well as transformation of the obtained solution into the solution of the source LPP are shown below. Picture 3 contains a detailed scheme of solving LPP using the game iteration method for the *modelling and applied problem solution* modes. Picture 4 contains a detailed algorithm showing LPP's solution scheme based on the game iteration method.

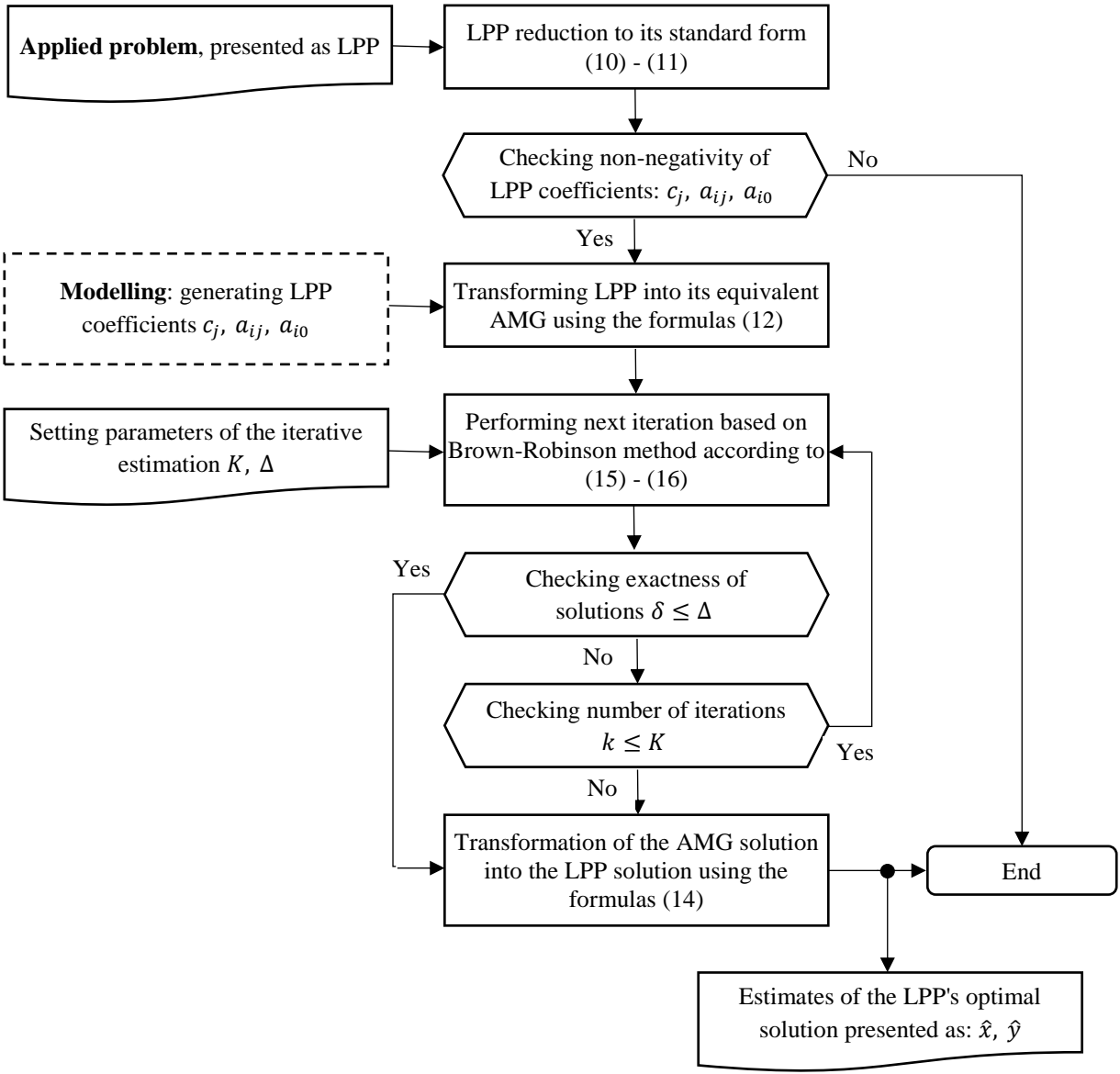

Picture 3 - Logic of the LPP algorithm solution scheme based on the game iteration method

---

**Algorithm 1.** LPP Solution algorithm using a game iteration method

---

**Given:** LPP at *max*: $L(\bar{y}) = \sum_{j=1}^{n} c_j y_j \to \max_{y_j}$ ; $\Omega_{\bar{y}} = \{\sum_{j=1}^{n} a_{ij} y_j \leq a_{i0};\ i = \overline{1,m};\ y_j \geq 0\}$.

All coefficients are set and positive: $c_j > 0,\ a_{ij} > 0,\ a_{i0} > 0$.

Iterative process stopping parameters are set: maximum permissible number of iterations is $K$; $\Delta$ is a threshold value of $\hat{\delta}_k$ measure.

**To find:** Optimal coordinate values of the vector of variables $\bar{y}$: $\bar{y}_{opt} = arg \max_{\Omega_{\bar{y}}} L(\bar{y})$.

- *Transforming LPP into AMG*

1: **for** $i \coloneqq 1$ **to** $m$ **do** {*Computing payoff matrix elements*}
2:    **for** $j \coloneqq 1$ **to** $n$ **do**
3:       $\tilde{a}_{ij} \coloneqq \frac{a_{ij}}{a_{i0} c_j}$
4:    **end for**
5: **end for**
6: $\tilde{A} \leftarrow [\tilde{a}_{ij}]_{mn}$ {*Formation of a payoff matrix*}

- *Solving AMG using a Brown-Robinson method*

7: $i(k) \coloneqq 1$ {*At the 1st iteration (k=1) related to P1, any pure strategy should be chosen, let it be the 1st*}
8: **for** $j \coloneqq 1$ **to** $n$ **do** {*A formation cycle of current average payoffs of P1*}

9:     $b_j \coloneqq \tilde{a}_{ij}$
10: $B_j(k) \coloneqq b_j$ {$B_j(k)$ is a current average accumulated payoff of P1 at k-numbered iteration}
11: **end for**
12: $V_{min}(k) \coloneqq \min_j B_j(k)$ {Current lower estimate of the game price}
13: $j(k) \coloneqq \arg\min_j B_j(k)$ {Number of P2's pure strategy is chosen according to the P1's payoff}
14: **for** $i \coloneqq 1$ **to** $m$ **do** {A formation cycle of current average payoffs of P2}
15:    $a_i \coloneqq \tilde{a}_{ij}$
16:    $A_i(k) \coloneqq a_i$ {$A_i(k)$ is a current average accumulated payoff of P2 at k-numbered iteration}
17: **end for**
18: $V_{max}(k) \coloneqq \max_i A_i(k)$ {Current upper estimate of the game price}
19: $V(k) \coloneqq (V_{min}(k) + V_{max}(k))/2$ {Current average estimate of the game price}
20: **for** $i \coloneqq 1$ **to** $m$ **do** {Formation of initial estimates of P1's mixed strategies}
21:    $\hat{p}_i \coloneqq 0$
22: **end for**
23: $\hat{p}_{i(k)} \coloneqq 1$
24: **for** $j \coloneqq 1$ **to** $n$ **do** {Formation of initial estimates of P2's mixed strategies}
25:    $\hat{q}_j \coloneqq 0$
26: **end for**
27: $\hat{q}_{j(k)} \coloneqq 1$
28: $k \coloneqq 2$ {Initial number of iteration related to the fictitious play cycle}
29: **do while** $\hat{\delta}_k > \Delta$ {AMG's fictitious play cycle until the set precision is achieved}
30:    $i(k) \coloneqq \arg\max_i A_i(k-1)$ {Number of P1's pure strategy is chosen according to the P2's payoff}
23:    **for** $j \coloneqq 1$ **to** $n$ **do** {A formation cycle of current average payoffs of P1}
24:       $b_j \coloneqq \tilde{a}_{ij}$
25:       $B_j(k) \coloneqq \left((k-1)B_j(k-1) + b_j\right)/k$ {Current average accumulated payoff of P1}
26:    **end for**
27:    $V_{min}(k) \coloneqq \min_j B_j(k)$ {Current lower estimate of the game price}
28:    $j(k) \coloneqq \arg\min_j B_j(k)$ {Number of P2's best strategy is chosen according to the P1's payoff}
29:    **for** $i \coloneqq 1$ **to** $m$ **do** {A formation cycle of current average payoffs of P2}
30:       $a_i \coloneqq \tilde{a}_{ij}$
31:       $A_i(k) \coloneqq \left((k-1)A_i(k-1) + a_i\right)/k$ {Current average accumulated payoff of P2}
32:    **end for**
33:    $V_{max}(k) \coloneqq \max_i A_i(k)$ {Current upper estimate of the game price}
34:    $V(k) \coloneqq (V_{min}(k) + V_{max}(k))/2$ {Current average estimate of the game price}
35:    **for** $i \coloneqq 1$ **to** $m$ **do** {Cycle of computing initial estimates of P1's mixed strategies}
36:       **if** $i = i(k)$ **then** $\hat{p}_{i(k)} \coloneqq \left((k-1)\hat{p}_{i(k-1)} + 1\right)/k$
37:          **else** $\hat{p}_{i(k)} \coloneqq \left((k-1)\hat{p}_{i(k-1)}\right)/k$
38:       **end if**
39:    **end for**
40:    **for** $j \coloneqq 1$ **to** $n$ **do** {Cycle of computing initial estimates of P2's mixed strategies}
41:       **if** $j = j(k)$ **then** $\hat{q}_{j(k)} \coloneqq \left((k-1)\hat{q}_{j(k-1)} + 1\right)/k$
42:          **else** $\hat{q}_{j(k)} \coloneqq \left((k-1)\hat{q}_{j(k-1)}\right)/k$
43:       **end if**
44:    **end for**
45:    $\hat{\delta}_k \coloneqq \max(\|\hat{p}_k - \hat{p}_{k-1}\|; \|\hat{q}_k - \hat{q}_{k-1}\|)$ {Convergence measure of the strategy estimates}
46:    **if** $k = K$ **then** {Checking whether the maximum permissible number of iterations was achieved}
47:       **exit do** {Exiting the fictitious play cycle on the basis of the iteration number}
48:    **end if**
49: **loop** {End of AMG's fictitious play cycle}
50: $\hat{p} \coloneqq \hat{p}_k$; $\hat{q} \coloneqq \hat{q}_k$; $\hat{V} \coloneqq V(k)$ {Final estimates of AMG solution}
- **Transformation of the AMG solution into the LPP solution**
51: $\hat{L}_{opt} \coloneqq \frac{1}{\hat{V}}$ {Estimate of an optimal value of the LPP's objective function}
52: **for** $i \coloneqq 1$ **to** $m$ **do** {Estimates of optimal values of variables related to the direct LPP}

53:     $\hat{x}_i^{opt} := \frac{\hat{p}_i}{a_{i0}\hat{v}}$
54: **end for**
55: **for** $j := 1$ **to** $n$ **do** {*Estimates of optimal values of variables related to the dual LPP*}
56:     $\hat{y}_j^{opt} := \frac{\hat{q}_j}{c_j\hat{v}}$
57: **end for**

Picture 4 - LPP algorithm solution scheme based on the game iteration method.

Algorithm 1 iterations are performed within a game block (lines 29-50). At that, the stopping rule is built according to the game indicators. However, if we move AMG-LPP solution transformation (lines 51-57) into the game block (after line 44), the convergence control can be performed directly according to the LPP solution estimates.

## COMPUTER IMPLEMENTATION AND TESTING

The game iteration method that was suggested in this work for LPP solution, was implemented as a separate Skat program, written in Object Pascal in Delphi environment. Picture 5 shows an interface window of the program in question (result of LPP solution).

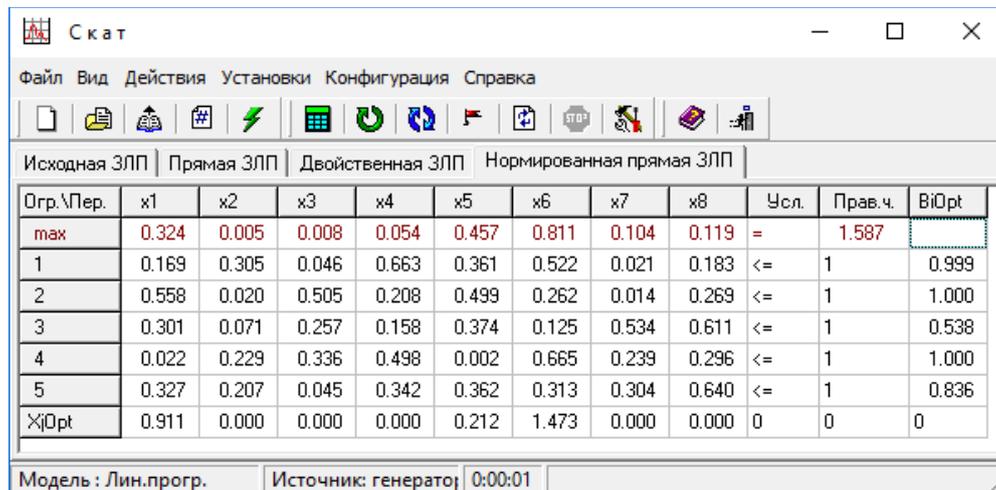

Pic. 5. LPP solution algorithm window using a game iteration method in Skat program

Skat allows solving rather large LPPs (up to several hundreds of variables and delimitations) with the necessary precision and within acceptable time (at modern computers - from fractions to unities of a second). Picture 5 shows results of solving LPP containing 12 variables and 15 delimitations. The program possesses large capabilities when downloading data from various sources and making necessary interface settings. When using two variables in LPP, Skat allows showing the results in a graphic form (see Pic. 6).

*Example.* Let us demonstrate a provided scheme of solving LPP using a game iteration method with the aid of Skat on a simple example of a limited resource distribution [2].

Let us presume the existence of a maximization task of some effect $L(\bar{y})$ that depends on the vector of variables $\bar{y} = [y_1 \quad y_2 \quad y_3]^T$. Structurally, these tasks correspond to the (7) - (8) setting, i.e. its optimal solution, in accordance with the game iteration method, shall be determined by the optimal strategy of a second player.

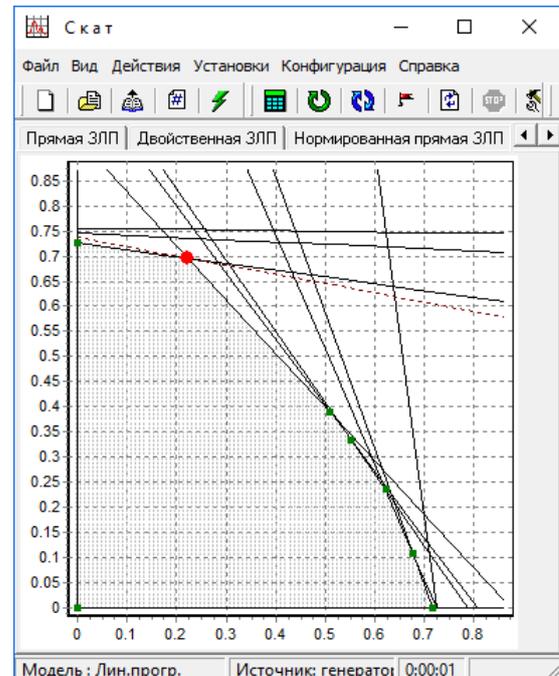

Pic. 6. Window with a graphic presentation of LPP solution

Therefore, according to the accepted definitions, the LPP is a dual one. Its objective function, criterion and delimitations are the following:

$$L(\bar{y}) = 3y_1 + 2y_2 + 5y_3 \to \max_{\bar{y}}, \quad (19)$$

$$\left.\begin{array}{l} y_1 + 2y_2 + y_3 \le 430 \\ 3y_1 + 0y_2 + 2y_3 \le 460 \\ y_1 + 4y_2 + 0y_3 \le 420 \end{array}\right\}, \quad (20)$$

$$y_j \ge 0, \ j = \overline{1,3}. \quad (21)$$

The antagonistic matrix game (19) - (21), being equivalent to LPP and with the consideration of (12), shall be presented as a following payoff matrix:

$$A = \begin{bmatrix} 77.519 & 232.558 & 46.512 \\ 217.391 & 0 & 86.957 \\ 79.365 & 476.190 & 0 \end{bmatrix} \times 10^{-5}. \quad (22)$$

Having solved the game with the aid of Brown-Robinson method in Skat, we receive the following:

$$\left.\begin{array}{l} \bar{p} = [0.218520 \ \ 0.681180 \ \ 0]^T \\ \bar{q} = [0 \ \ 0.148149 \ \ 0.851851]^T \\ V = 74.074 \times 10^{-5} \end{array}\right\} \quad (23)$$

The optimal mixed strategies of the second player ($\bar{q}$) correspond to the optimal solution of a source dual LPP, while the $\bar{p}$ strategy of the first player corresponds to the direct LPP.

Therefore, based on (14) and (23), the solution of a dual problem shall be the following:

$$y_1 = 0; \ y_2 = 100; \ y_3 = 230; \ L = 1350, \quad (24)$$

solution of a direct LPP:

$$x_1 = 1; \ x_2 = 2; \ x_3 = 0; \ W = 1350. \quad (25)$$

In order to check the appropriateness of the solution that was obtained by the game iteration method, the optimal solution of the (19) - (21) problem was also found with the aid of Solution Finder add-in in MS Excel. The obtained solution of a direct and dual LPP was absolutely the same as the values (24) and (25).

## RESULTS AND DISCUSSION

1. The main characteristic of a *game iteration method* is a high rate of the optimal solution receipt. According to the data provided by the bibliographical sources [3, 8, 12-14], simplex method allows looking for the solution with the rate of the variables that depend exponentially on the space dimensionality (for instance, like $c \times 2^n$ where $c$ is an invariable and $n$ is a space dimensionality of the variables), while Karmarkar method provides a polynomial dependency of the solution's search rate on the space dimensionality (for example, $c \times n \times m$).

The simulation experiment examining dependency of the convergence rate related to the LPP solution estimates, on the problem parameters has been performed with the aid of Skat software using a modern computer with the following parameters: Intel Core i7-4702MQ CPU 2*2.4 GHz, RAM 16 Gb. The analysis performed with a fixed number of iterations ($10^6$) demonstrated that the computational speed (time $t$) does not depend only on the number of LPP variables ($n$) and number of delimitations ($m$), being first of all dependent on their total ($m + n$). The regression dependency is the following: $t = 0.044(m + n)$; [sec]. For example, the computational time for $m + n = 100 + 100 = 200$ is about 8.8 seconds. It is necessary to mention that number of iterations equal to $10^6$ provides a very high computational precision even for large problems. The precision in question is usually acceptable for the majority of applied problems on optimal control. Thus, when talking about the above-mentioned example, the precision constitutes fractions and units of percents even with the iteration rate equal to $10^2$.

If the system response time is important, it is desirable to know the dependency of computational time on the number of iterations. We have varied this parameter in the computing experiment, where a degree ($e$) of 10 was considered as a factor. The dimensionality factor ($m + n$) has also been varied. Quadratic (based on $e$ factor) two-factor regression model can be presented as: $t = 0.037(m + n) - 3.82e + 0.56e^2$; [sec]. The determination coefficient of the model is $R^2 = 0.93$. It provides a good description of the dependency between the computing duration and number of iterations within the $10^5 - 10^7$ range. With the values less than $10^5$, even the large dimensionality computation is performed within deciseconds or centiseconds. Certainly, when implementing the game iteration algorithm on some on-board computers or microprocessors it is necessary to perform additional research devoted to the dependency of the computing time on the parameters of the computing environment and the problems that are being solved.

2. The advantages of the *game iteration method* include its simplicity, transparency and a relatively high rate of the solution receipt. We can control the solution receipt rate by choosing an acceptable level of the calculation precision. These parameters open new possibilities for using linear optimization in the integrated control systems that are widely applied nowadays, for instance in robotic systems, smart homes, Internet of things etc.

3. Disadvantages of the *game iteration method* include some limitation of the solved problems, caused by the fact that they should possess non-negative values of the right parts of delimitations and objective function coefficients. In some cases, this delimitation can be weakened by linear transformation of a variables space. However, it is necessary to mention that a significant number of problems on the optimal distribution of limited resources are related to the specific type that is discussed in this work.

## CONCLUSIONS

1. The game iteration method that has been suggested in this work and oriented at the computational solution of a widely applied in practice optimization problem, i.e. linear programming problem, proved its work capacity with the problems of large dimensionality.

2. Using multiple model examples, the author has showed that the rate of the solution's computational search with the application of the game iteration method is linearly dependent on the total amount of dimensionality of the variables space and number of delimitations, which is an important advantage of the suggested method in comparison with others that are used when solving linear programming problems.

3. Implementation of a game iteration method in the Skat program that is presented in the work, allows solving linear programming problems in various applications, while its libraring version (*.dll*) provides a possibility for other applications to import functions of a game iteration method.

4. High rate of solving linear programming optimization problems using a game iteration method allows using this easy and stable algorithm as a part of on-board algorithms responsible for optimal management of modern integrated control systems in various applied spheres.

## REFERENCES

[1] Mine H., Osaki S. Markovian Decision Processes. Moscow: Nauka Publishing House. 1977. p. 176
[2] V. Ya. Vilisov. Adaptive Choice of Managerial Decisions. Operation Examination Models as the Means of Storing the Decision Taker's Knowledge - Saarbruecken (Germany): LAP LAMBERT Academic Publishing, 2011, p. 376.
[3] Taha H.A. Operation research: An Introduction. 8-th ed., Pearson Prentice Hall, New Jersey, 2007, p. 813
[4] John von Neumann and Oskar Morgenstern, Theory of Games and Economic Behavior. Moscow: Nauka Publishing House, 1970. p. 707
[5] Guillermo Owen. Game Theory - Academic Press, 1982, p. 344
[6] On Robot Algorithms Adapted to the Target Preferences of the Decision Taker // Book of reports of the Russian National Scientific and Technological Conference "Extreme Robotics", Saint Petersburg: Politekhnika-Service, 2012, pp. 124-130
[7] V.I. Danilin, Financial and Operation Planning in the Company. Methods and Models: textbook. Moscow: Delo Publishing House of RANEPA. 2014. p. 616
[8] George Bernard Dantzig, Linear Programming: Theory and Extensions. Moscow: Progress Publishing House, 1966. p. 601
[9] Brown G.W. Iterative solution of games by fictitious play. In Activity Analysis of Production and Allocation, Cowles Commission Monograph No. 13, pp. 374–376. Wiley, New York, N.Y. (1951)
[10] Robinson J. An iterative method of solving a game. Ann. Math. (2) 54, 296–301 (1951).
[11] O.O. Iemets, D.M. Olhovskiy, O.V. Olhovskaya, Comparison of Methods for Solving Game Tasks: Numerical Experiments // Artificial Intelligence. 2014. No. 1. pp. 47-56. URL: http://dspace.nbuv.gov.ua/bitstream/handle/123456789/85236/7-Iemets.pdf?sequence=1.
[12] Kobzar A.I., Tikmenova I.V., Tikmenov V.N., Comparative Analysis of m×n Matrix Game Solution by Linear Programming Method and Brown - Robinson Iteration Method // Electronic Information Systems. 2014. No. 3 (3). pp. 33-52. URL: https://elibrary.ru/item.asp?id=23618071.
[13] Emrah A., Handan A., Serkan D. Brown–Robinson method for interval matrix games // Soft Comput. 2011. No.15. pp. 2057-2064. DOI 10.1007/s00500-011-0703-6.
[14] Akyar H. Fuzzy Risk Analysis for a Production System Based on the Nagel Point of a Tri-angle, Mathematical Problems in Engineering, vol. 2016, pp. 1–9, 2016.